\documentclass[preprint,11pt]{elsarticle}
\usepackage{graphicx,graphics}
\usepackage[intlimits]{amsmath}
\usepackage{amssymb}
\usepackage{exscale}
\usepackage{times}
\usepackage{subfigure}
\usepackage{subfig}
\usepackage{multirow}
\usepackage[abs]{overpic}

\usepackage{pgfplots}
\usetikzlibrary{plotmarks}

\usepackage{geometry}

\journal{Composite Structures}

\newcommand{\Eref}[1]{Equation (\ref{#1})}
\newcommand{\fref}[1]{Figure (\ref{#1})}
\newcommand{\Erefs}[1]{Equations (\ref{#1})}
\newcommand{\frefs}[1]{Figures~(\ref{#1})}

\newcommand{\rmd}{\rm{d}}

\newcommand{\T}[0]{\mathrm{T}}

\begin{document}

\begin{frontmatter}

\title{Free vibration and mechanical buckling of plates with in-plane material inhomogeneity - a three dimensional consistent approach}      

\author[unsw]{Tingsong Xiang}
\author[unsw]{Sundararajan Natarajan \corref{cor1}\fnref{fn1}}
\author[unsw]{Hou Man}
\author[unsw]{Chongmin Song}
\author[unsw]{Wei Gao}

\cortext[cor1]{Corresponding author}

\address[unsw]{School of Civil \& Environmental Engineering, The University of New South Wales, Sydney, Australia}

\fntext[fn1]{School of Civil and Environmental Engineering, University of New South Wales, Sydney, NSW 2052, Australia. Tel: +61 2 93855030, Email: snatarajan@cardiffalumni.org.uk; sundararajan.natarajan@gmail.com}

\begin{abstract}
In this article, we study the free vibration and the mechanical buckling of plates using a three dimensional consistent approach based on the scaled boundary finite element method.  The in-plane dimensions of the plate are modeled by two-dimensional higher order spectral element. The solution through the thickness is expressed analytically with Pad\'e expansion. The stiffness matrix is derived directly from the three dimensional solutions and by employing the spectral element, a diagonal mass matrix is obtained. The formulation does not require ad hoc shear correction factors  and no numerical locking arises. The material properties are assumed to be temperature independent and graded only in the in-plane direction by a simple power law. The effective material properties are estimated using the rule of mixtures. The influence of the material gradient index, the boundary conditions and the geometry of the plate on the fundamental frequencies and critical buckling load are numerically investigated. 
\end{abstract}


\end{frontmatter}

\section{Introduction}
\vspace{-6pt}

The introduction of new class of engineered materials, coined as functionally graded materials (FGMs) has spurred the interest among researchers to study the response of structures with these materials. The FGMs are characterized by \emph{smooth and continuous} transition of material properties from one surface to another. Typically FGMs are made from a mixture of a ceramic and metal. The ceramic constituent provides thermal stability due to its low thermal conductivity, whilst the metallic phase provides structural stability. FGMs eliminate the sharp interfaces existing in laminated composites with a gradient interface and are considered to be an alternative material in many engineering applications. The material properties can be graded in the thickness direction, in the in-plane or in both the directions. It can be seen from the literature that considerable attention has been devoted to functionally graded material plates with properties graded in the thickness direction~\cite{reddy2000,ganapathiprakash2006,ferreirabatra2006,reddychin1998}. From the literature, it can be seen that the static and the dynamic response of functionally graded material plates and shells is studied extensively. It is beyond the scope of this paper to review the literature on plate/shells with material properties graded in the thickness direction. Interested readers are referred to the literature and references therein and a recent review by Jha and Kant~\cite{jhakant2013}. To the author's knowledge there are only a few investigations on the structural response of structures in which the material is graded in the in-plane direction or in both directions~\cite{nemat-alla2003,qianching2004,qianbatra2005,goupeevel2006,luchen2008,liuwang2010,uymazaydogdu2012}. 
Qian and Ching~\cite{qianching2004} and Qian and Batra~\cite{qianbatra2005} optimized the fundamental frequency of bi-directional\footnote{material properties graded in the thickness and in the in-plane direction} functionally graded beams and plates by employing meshless local Petrov Galerkin method. By employing element free Galerkin method, Goupee and Vel~\cite{goupeevel2006} optimized the natural frequency of bidirectional functionally graded beams. Nemat-Alla~\cite{nemat-alla2003} by employing the rule of mixtures studied the thermal response of FGM structures graded in both the directions. L\"u \textit{et al.,}~\cite{luchen2008} derived semi-analytical solutions based on differential quadrature method for beams graded in both the directions. It was observed that the thermal stresses can be reduced by bi-directional functionally gradation instead of the conventional unidirectional functionally graded materials. Very recently, Liu \textit{et al.,}~\cite{liuwang2010} and Uymaz \textit{et al.,}~\cite{uymazaydogdu2012} studied the fundamental frequency of plates with in-plane material inhomogeneity by Levy's type solution and differential quadrature method, respectively. It was observed that the fundamental frequency of the plate with in-plane material inhomogeneity is highly influenced by the gradation.

It is noted that in all of the above studies, the plate structures are modelled by employing two dimensional structural theories. The different approaches employed are: single layer theories, discrete layer theories and mixed plate theory. In the single layer theory approach, the plate is assumed to be one equivalent single layer (ESL), whereas in the discrete layer theory, each layer is considered, for example in the case of laminated composites. Although the discrete layer theories provide very accurate results, increasing the number of layers increases the number of unknowns and in turn the computational time. Recently, Carrera~\cite{carrerademasi2002,carrera2003} derived a series of axiomatic approaches for general description of two-dimensional formulations for multilayered plates and shells. With this unified formulation, it is possible to implement in a single software a series of hierarchical formulations, thus affording a systematic assessment of different theories ranging from simple ESL models up to higher order layerwise descriptions. The aforementioned plate theories have been used to develop discrete models such as the finite element method~\cite{greimannlynn1970,bathedvorkin1985,nguyen-xuantran2012}, meshless methods~\cite{ferreirabatra2005,ferreirabatra2006,liewzhao2011} and more recently, isogeometric analysis~\cite{thainguyen-xuan2012,valizadehnatarajan2013,tranferreira2013}. A comprehensive review of various meshless methods for analyzing plates and shells is given in~\cite{liewzhao2011}. Although, the numerical methods provide a general and systematic technique to analyze plate structures, difficulties still exist in the development of plate elements based on the above mentioned plate theories, one of which is the shear locking phenomenon. It can be seen that considerable effort has been devoted to suppress shear locking~\cite{bathedvorkin1985,brezzibathe1989,hughescohen1978,somashekarprathap1987}.

\paragraph{Approach} However, plates are essentially three dimensional structures. For predicting the realistic behaviour, more accurate analytical/numerical models based on the three-dimensional models are required. In this paper, we study the free vibration and mechanical buckling of plates with in-plane inhomogeneity using a recently developed three-dimensional consistent approach \cite{mansong2012}. This approach is based on the scaled boundary finite element method \cite{songwolf1997a}. 
The formulations are directly derived from three-dimensional governing equations without \emph{any plate assumptions}. Only the in-plane dimensions of the plate are discretised and any displacement-based elements can be used. The stiffness matrix is derived from the three dimensional solution, which is expressed analytically in the through-thickness direction with Pad\'e expansion. Thus, no numerical locking arises. The use of high-order spectral elements leads to an efficient stiffness matrix construction. A diagonal mass matrix is also derived such that the free vibration in our study is expressed as a standard eigenproblem.

\paragraph{Outline} The paper is organized as follows. In the following section, after discussing the functionally graded material, the three dimensional consistent approach to analyse plate structures is presented. Section \ref{elemdesc} describes the element employed in this study. The numerical results for the free vibration and critical buckling of thin functionally graded material plates are given in Section \ref{numres}, followed by concluding remarks in the last section.

\section{Theoretical Formulation}
\label{theorydev}

\subsection{Functionally graded material} 
\begin{figure}[htpb]
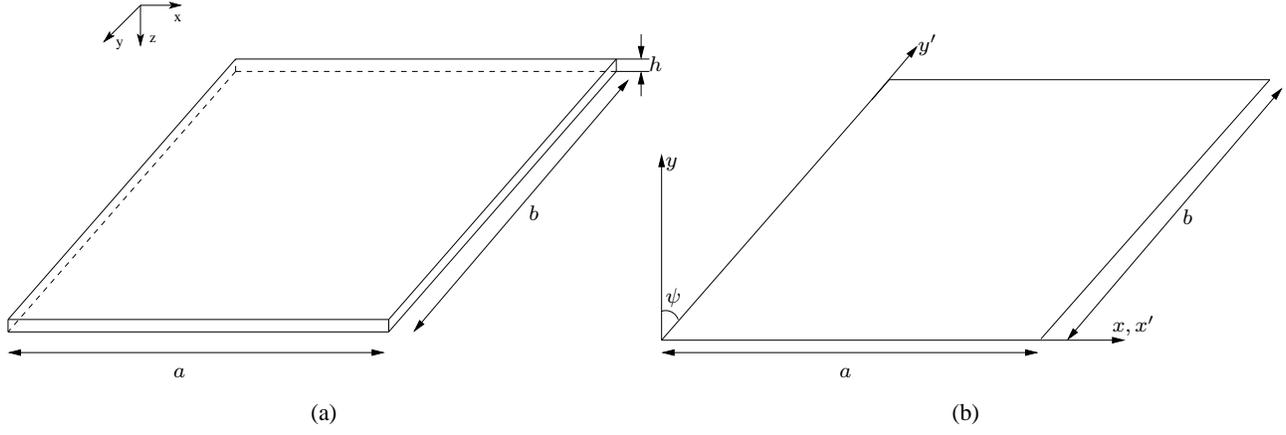

\centering
\subfigure[]{\input{./Figures/plate.pstex_t}}
\subfigure[]{\input{./Figures/skew.pstex_t}}
\caption{Coordinate system of a rectangular skew plate.}
\label{fig:platefig}
\end{figure}

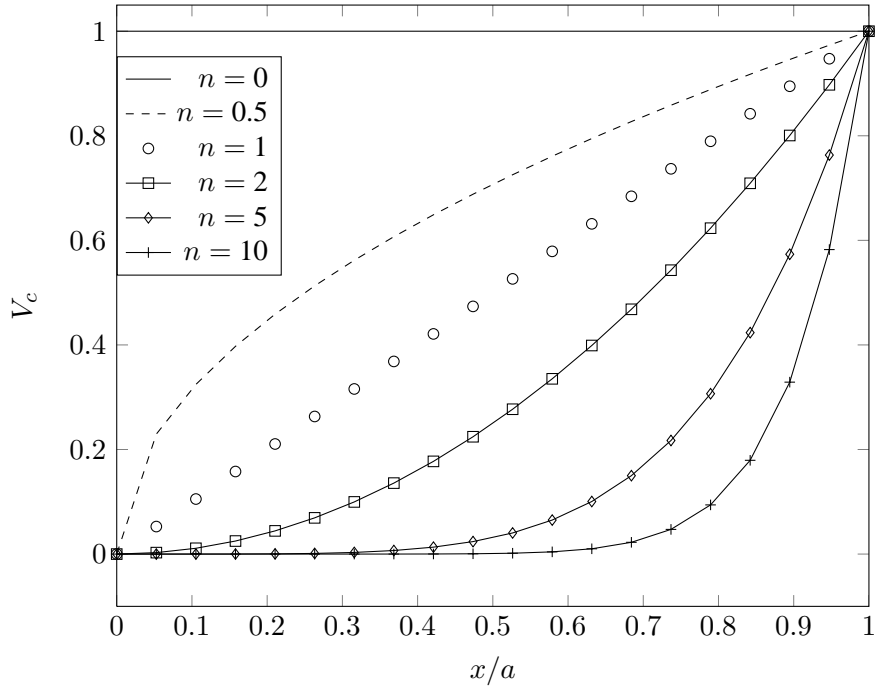
\begin{figure}[htpb]
\centering
\newlength\figureheight 
\newlength\figurewidth 
\setlength\figureheight{8cm} 
\setlength\figurewidth{10cm}
%
%
%
%
\begin{tikzpicture}

\begin{axis}[%
width=\figurewidth,
height=\figureheight,
scale only axis,
xmin=0,
xmax=1,
xlabel={$x/a$},
ymin=-0.1,
ymax=1.05,
ylabel={$V_c$},
legend style={at={(0,0.55)},anchor=south west,draw=black,fill=white,legend cell align=right}
]
\addplot [
color=black,
solid
]
table[row sep=crcr]{
0 1\\
0.0526315789473684 1\\
0.105263157894737 1\\
0.157894736842105 1\\
0.210526315789474 1\\
0.263157894736842 1\\
0.315789473684211 1\\
0.368421052631579 1\\
0.421052631578947 1\\
0.473684210526316 1\\
0.526315789473684 1\\
0.578947368421053 1\\
0.631578947368421 1\\
0.684210526315789 1\\
0.736842105263158 1\\
0.789473684210526 1\\
0.842105263157895 1\\
0.894736842105263 1\\
0.947368421052632 1\\
1 1\\
};
\addlegendentry{$n=$ 0};

\addplot [
color=black,
dashed
]
table[row sep=crcr]{
0 0\\
0.0526315789473684 0.229415733870562\\
0.105263157894737 0.324442842261525\\
0.157894736842105 0.397359707119513\\
0.210526315789474 0.458831467741124\\
0.263157894736842 0.512989176042577\\
0.315789473684211 0.561951486949016\\
0.368421052631579 0.606976978666884\\
0.421052631578947 0.64888568452305\\
0.473684210526316 0.688247201611685\\
0.526315789473684 0.725476250110012\\
0.578947368421053 0.760885910252682\\
0.631578947368421 0.794719414239026\\
0.684210526315789 0.827170191868511\\
0.736842105263158 0.858395075278952\\
0.789473684210526 0.888523316638639\\
0.842105263157895 0.917662935482247\\
0.894736842105263 0.945905302926917\\
0.947368421052632 0.973328526784575\\
1 1\\
};
\addlegendentry{$n=$ 0.5};

\addplot [
color=black,
only marks,
mark=o,
mark options={solid}
]
table[row sep=crcr]{
0 0\\
0.0526315789473684 0.0526315789473684\\
0.105263157894737 0.105263157894737\\
0.157894736842105 0.157894736842105\\
0.210526315789474 0.210526315789474\\
0.263157894736842 0.263157894736842\\
0.315789473684211 0.315789473684211\\
0.368421052631579 0.368421052631579\\
0.421052631578947 0.421052631578947\\
0.473684210526316 0.473684210526316\\
0.526315789473684 0.526315789473684\\
0.578947368421053 0.578947368421053\\
0.631578947368421 0.631578947368421\\
0.684210526315789 0.684210526315789\\
0.736842105263158 0.736842105263158\\
0.789473684210526 0.789473684210526\\
0.842105263157895 0.842105263157895\\
0.894736842105263 0.894736842105263\\
0.947368421052632 0.947368421052632\\
1 1\\
};
\addlegendentry{$n=$ 1};

\addplot [
color=black,
solid,
mark=square,
mark options={solid}
]
table[row sep=crcr]{
0 0\\
0.0526315789473684 0.00277008310249307\\
0.105263157894737 0.0110803324099723\\
0.157894736842105 0.0249307479224377\\
0.210526315789474 0.0443213296398892\\
0.263157894736842 0.0692520775623269\\
0.315789473684211 0.0997229916897507\\
0.368421052631579 0.135734072022161\\
0.421052631578947 0.177285318559557\\
0.473684210526316 0.224376731301939\\
0.526315789473684 0.277008310249307\\
0.578947368421053 0.335180055401662\\
0.631578947368421 0.398891966759003\\
0.684210526315789 0.46814404432133\\
0.736842105263158 0.542936288088643\\
0.789473684210526 0.623268698060942\\
0.842105263157895 0.709141274238227\\
0.894736842105263 0.800554016620499\\
0.947368421052632 0.897506925207756\\
1 1\\
};
\addlegendentry{$n=$ 2};

\addplot [
color=black,
solid,
mark=diamond,
mark options={solid}
]
table[row sep=crcr]{
0 0\\
0.0526315789473684 4.03861073406192e-07\\
0.105263157894737 1.29235543489982e-05\\
0.157894736842105 9.81382408377048e-05\\
0.210526315789474 0.000413553739167941\\
0.263157894736842 0.00126206585439435\\
0.315789473684211 0.00314042370680655\\
0.368421052631579 0.00678769306073788\\
0.421052631578947 0.0132337196533741\\
0.473684210526316 0.0238475925235623\\
0.526315789473684 0.0403861073406192\\
0.578947368421053 0.0650422297331407\\
0.631578947368421 0.10049355861781\\
0.684210526315789 0.149950789528205\\
0.736842105263158 0.217206177943612\\
0.789473684210526 0.306682002617828\\
0.842105263157895 0.423479028907972\\
0.894736842105263 0.573424972103296\\
0.947368421052632 0.763122960753992\\
1 1\\
};
\addlegendentry{$n=$ 5};

\addplot [
color=black,
solid,
mark=+,
mark options={solid}
]
table[row sep=crcr]{
0 0\\
0.0526315789473684 1.63103766612802e-13\\
0.105263157894737 1.67018257011509e-10\\
0.157894736842105 9.63111431471934e-09\\
0.210526315789474 1.71026695179785e-07\\
0.263157894736842 1.59281022082814e-06\\
0.315789473684211 9.86226105827261e-06\\
0.368421052631579 4.60727770867891e-05\\
0.421052631578947 0.0001751313358641\\
0.473684210526316 0.000568707669169863\\
0.526315789473684 0.00163103766612802\\
0.578947368421053 0.00423049164865866\\
0.631578947368421 0.0100989553236712\\
0.684210526315789 0.0224852392801322\\
0.736842105263158 0.0471785237368721\\
0.789473684210526 0.0940538507296812\\
0.842105263157895 0.179334487924839\\
0.894736842105263 0.328816198631666\\
0.947368421052632 0.582356653229939\\
1 1\\
};
\addlegendentry{$n=$ 10};

\end{axis}
\end{tikzpicture}%
\caption{Volume fraction of the ceramic phase as a function of global $x$-coordinate.}
\label{fig:ceramicvolfrac}
\end{figure}

Consider a functionally graded material (FGM) rectangular skew plate with length $a$, width $b$, height $h$ and skew angle $\psi$  made by mixing two distinct material phases, viz., a ceramic phase and a metallic phase. The ceramic phase provides thermal stability, whilst the metallic phase provides structural stability. Assume the coordinates $x,y$ along the in-plane directions and $z$ along the thickness direction (see \fref{fig:platefig}). The material is assumed to be graded only in the in-plane direction (along global $x$) according to a power law distribution whilst it is constant through the thickness direction. The homogenized material properties can be computed by employing the rule of mixtures.  The effective Young's modulus $E$, Poisson's ratio $\nu$ and the mass density $\rho$ of the FGM, evaluated using the rule of mixtures are:
\begin{align}
E &= V_m E_m + V_c E_c \nonumber \\
\nu &= V_m \nu_m + V_c \nu_c \nonumber \\
\rho &= V_m \rho_m + V_c \rho_c
\label{fgmc}
\end{align}
Here $V_i (i=c,m)$ is the volume fraction of the phase material. The subscripts $c$ and $m$ refer to the ceramic and metal phases, respectively. The volume fractions are related by $V_c + V_m = 1$ and $V_c$ is expressed as
\begin{equation}
V_c = (x/a)^n
\label{eqn:volumefraccer}
\end{equation}
where $n$ in \Eref{eqn:volumefraccer} is the volume fraction exponent, also referred to as the material gradient index in the literature, and $x$ is referred to the global coordinate system. \fref{fig:ceramicvolfrac} shows the variation of the volume fraction of ceramic phase along the in-plane directions.

\vspace{10pt}
\subsection{3D governing equations for plate structures}
Consider a plate of constant thickness $h$ and with length $a$ and width $b$ (see \fref{fig:platefig}). The displacement components along the $(x,y)$ directions and $z-$direction are denoted as $u_{x}=u_{x}(x,y,z)$, $u_{y}=u_{y}(x,y,z)$ and $u_{z}=u_{z}(x,y,z)$. The displacement vector $\mathbf{u}=\mathbf{u}(x,y,z)$ is arranged as $\mathbf{u}=[u_{z},u_{x},u_{y}]^{\rm T}$. The strains $\{\varepsilon\}=\{\varepsilon(x,\, y,\, z)\}$
are expressed as 
\begin{align}
\boldsymbol{\varepsilon} & =[\varepsilon_{z}\;\;\varepsilon_{x}\;\;\varepsilon_{y}\;\;\gamma_{xy}\;\;\gamma_{yz}\;\;\gamma_{xz}]^{\T}=\mathbf{L} \mathbf{u},\label{straineq-1}
\end{align}
where $\mathbf{L}$ is the differential operator. The stresses $\boldsymbol{\sigma} =\{\sigma(x,\, y,\, z)\}$ follow from Hooke's
law with the elasticity matrix $\mathbf{D}$ as 
\begin{equation}
\boldsymbol{\sigma} = [\sigma_{z}\;\;\sigma_{x}\;\;\sigma_{y}\;\;\tau_{xy}\;\;\tau_{yz}\;\;\tau_{xz}]^{\T}=\mathbf{D} \boldsymbol{\varepsilon}.
\label{eq-stress_strain}
\end{equation}
The equation of equilibrium with vanishing body force is written as
\begin{equation}
\mathbf{L}^{\T} \boldsymbol{\sigma}  = \rho\ddot{\mathbf{u}} \label{eq:equilibrium}
\end{equation}
and the bottom and top surfaces of the plate may be subjected to surface traction. The strain energy $U$, the work done by the applied external forces $V$ and the kinetic energy $T$ is given by:
\begin{align}
U &= \int\limits_\Omega  \left[ \varepsilon_{x} \sigma_{x} + \varepsilon_{y} \sigma_{y} + \gamma_{xy} \tau_{xy} + \gamma_{xz} \tau_{xz} + \gamma_{yz} \tau_{yz} + \varepsilon_{z} \sigma_{z} \right]  ~\rmd \Omega \nonumber \\
V &= \int\limits_\Omega \left[ N_{x} (u_{z,x})^2 + N_{y}  (u_{z,y})^2 + 2N_{xy} u_{z,x}u_{z,y} \right] ~\rmd \Omega \nonumber \\
T &= \int\limits_\Omega \left. \rho[ \delta(u_{x})\ddot{u}_{x}+ \delta(u_{y})\ddot{u}_{y} + \delta(u_{z})\ddot{u}_{z} \right]~\rmd \Omega
\end{align}
The derivation of the governing equations is based on the principle of virtual work equation
\begin{equation}
\delta( U + V - T) = 0
\end{equation}
In the following section, the scaled boundary finite element method will be employed to derive the stiffness matrix for the plate structure. The conventional finite element procedure~\cite{zinekiewicztaylor2000} is employed to derive the mass matrix $\mathbf{M}$ and the geometric stiffness matrices, $\mathbf{K}_G$. 

\subsection{3D consistent approach for plate structures} 
The geometry of the plate is modeled by translating the two dimensional (2D) mesh along the $z-$direction, where the geometry of an in-plane 2D plate elements is obtained by interpolating the nodal coordinates $\mathbf{x}$ and $\mathbf{y}$ using the shape function $\mathbf{N}(\eta,\zeta)$ formulated in the local coordinate $\eta$ and $\zeta$
\begin{align}
\bf{x} (\eta,\zeta) & =\mathbf{N} (\eta,\zeta) \mathbf{x}
\end{align}
In this study, $\mathbf{N}$ is based on high-order spectral elements that are detailed in the next section. However, other shape functions, such as the moving least square approximations (MLS) and non-uniform rational B-splines can be also employed to discretize the in-plane dimensions. 

The strain in Equation \eqref{straineq-1} is rewritten as
\begin{equation}
\boldsymbol{\varepsilon} = \mathbf{b}_1 \mathbf{u}_{,z} + \mathbf{b}_2 \mathbf{u}_{,\eta} + \mathbf{b}_3 \mathbf{u}_{,\zeta}
\end{equation} 
where
\begin{align}
\mathbf{b}_1 &=  \frac{1}{|J|} \left[\begin{array}{ccc}
1 & 0 & 0\\
0 & 0 & 0\\
0 & 0 & 0\\
0 & 0 & 0\\
0 & 0 & 1\\
0 & 1 & 0
\end{array}\right]; \nonumber \\
\mathbf{b}_2 &= \frac{1}{|J|}\left[\begin{array}{ccc}
0 & 0 & 0\\
0 & y_{,\zeta} & 0\\
0 & 0 & -x_{,\zeta}\\
0 & -x_{,\zeta} & y_{,\zeta}\\
-x_{,\zeta} & 0 & 0\\
y_{,\zeta} & 0 & 0
\end{array}\right]; \nonumber \\
 \mathbf{b}_3 &= \frac{1}{|J|}\left[\begin{array}{ccc}
0 & 0 & 0\\
0 & -y_{,\eta} & 0\\
0 & 0 & x_{,\eta}\\
0 & x_{,\eta} & -y_{,\eta}\\
x_{,\eta} & 0 & 0\\
-y_{,\eta} & 0 & 0
\end{array}\right].
\end{align}

The 3D displacement field of the plate, $\mathbf{u}$, is represented semi-analytically here. The displacement variations along the line pass through a node of the 2D mesh and normal to the mid-plane are expressed analytically by functions $\mathbf{u}(z)$ of the coordinate $z$. The 3D displacement field is described by interpolating the displacement functions $\mathbf{u}(z)$ using the same shape function $\mathbf{N}\equiv\mathbf{N}(\eta,\zeta)$ such that $\mathbf{u}=\mathbf{N}\mathbf{u}(z)$. By employing the principle of virtual work as shown in \cite{mansong2012} for the full derivation, the internal nodal force is derived into
\begin{equation}
\mathbf{q}(z)=\mathbf{E}_0 \mathbf{u_{\mathit{,z}}}(z)+\mathbf{E}_{1}^{\rm T} \mathbf{u}(z)
\label{eq:q1}
\end{equation}
Satisfying the virtual work equation also brings
\begin{equation}
\mathbf{q}_{,z}(z)=\mathbf{E}_1 \mathbf{u_{\mathit{,z}}}(z)+\mathbf{E}_2 \mathbf{u}(z)\label{eq:q2}
\end{equation}
where $\mathbf{E}_0$, $\mathbf{E}_1$ and $\mathbf{E}_2$ are the scaled boundary finite element coefficient matrices:
\begin{align}
\mathbf{E}_0 & =\int_{-1}^{+1}\int_{-1}^{+1}\mathbf{B}_{1}^{\rm T} \mathbf{D} \mathbf{B}_1 |J| ~\rmd \eta ~\rmd \zeta \nonumber \\
\mathbf{E}_1 & =\int_{-1}^{+1}\int_{-1}^{+1}\mathbf{B}_{2}^{\rm T} \mathbf{D} \mathbf{B}_1 |J| ~\rmd \eta ~\rmd \zeta \nonumber \\
\mathbf{E}_2 & =\int_{-1}^{+1}\int_{-1}^{+1}\mathbf{B}_{2}^{\rm T} \mathbf{D} \mathbf{B}_2 |J| ~\rmd \eta ~\rmd \zeta
\end{align}
with
\begin{equation}
\begin{array}[t]{cc}
\mathbf{B}_1=\mathbf{b}_1 \mathbf{N}; & \mathbf{B}_2 = \mathbf{b}_2 \mathbf{N_{\mathit{,\eta}}}+\mathbf{b}_3 \mathbf{N_{\mathit{,\zeta}}}\end{array}
\end{equation}
where $\mathbf{D}$ is the elasticity matrix 
\begin{equation}
\mathbf{D}=\frac{E(x)}{1-\nu^2(x)}\left[\begin{array}{ccc}
1 & \nu(x) & 0\\
\nu(x) & 1 & 0\\
0 & 0 & \frac{1-\nu(x)}{2}
\end{array}\right],
\end{equation}
and is evaluated with Equation~\eqref{fgmc}.
The determinant of the Jacobin matrix is given by $\mathbf{|\mathit{J}|}=x_{,\eta}y_{,\zeta}-x_{,\zeta}y_{,\eta}$. By introducing the variable
\begin{equation}
\mathbf{X}(z)=\left\{ \begin{array}{c}
\mathbf{u}(z)\\
\mathbf{q}(z)
\end{array}\right\} ,
\end{equation}
Equations~\eqref{eq:q1} and \eqref{eq:q2} are combined into  
\begin{equation}
\mathbf{X_{\mathit{,z}}}(z)=-\mathbf{A}\mathbf{X}(z)\label{eq:1stode}
\end{equation}
with the coefficient matrix
\begin{equation}
\mathbf{A}=\left[\begin{array}{cc}
\mathbf{A_{11}} & \mathbf{A_{12}}\\
\mathbf{A_{21}} & \mathbf{A_{22}}
\end{array}\right]=\left[\begin{array}{cc}
\mathbf{E}_{\mathbf{0}}^{-1}\mathbf{E}_{\mathbf{1}}^{\mathrm{T}} & -\mathbf{E}_{\mathbf{0}}^{-1}\\
-\mathbf{E_{2}}+\mathbf{E}_{\mathbf{1}}\mathbf{E}_{0}^{-1}\mathbf{E}_{\mathbf{1}}^{\mathrm{T}} & -\mathbf{E}_{\mathbf{1}}\mathbf{E}_{\mathbf{0}}^{-1}
\end{array}\right].\label{eq:Z}
\end{equation}
The general solution of $\mathbf{X}(z)$ is given as
\begin{equation}
\mathbf{X}(z)=\mathrm{e}^{-\mathbf{A}z}\mathbf{c}\label{fsol}
\end{equation}
By applying a Pad\'e expansion of order $(2,2)$ to express $\mathbf{X}(z)$ and substituting the boundary conditions at the top and bottom surfaces of the plate, a 3D consistent stiffness matrix is obtained \cite{mansong2013}. When further transforming the 3D displacements into typical plate degree of freedoms $\mathbf{d}=\left[\boldsymbol{\theta},\boldsymbol{u}\right]^{\mathrm{T}}$, the stiffness matrix for the plate is devised: 
\begin{equation}
\mathbf{K}=h\left[\begin{array}{cc}
\mathbf{E}_{\mathbf{0}}\left(\mathbf{I}+h^{2}\mathbf{V_{11}}\right) & \mathbf{E}_{\mathbf{1}}^{\mathrm{T}}\\
{}\mathbf{E}_{\mathbf{1}}\left(\mathbf{I}+h^{2}\mathbf{V_{11}}\right)-h^{2}\mathbf{V_{21}} & \mathbf{E}_{\mathbf{2}}
\end{array}\right]
\label{eq:fullplateK}
\end{equation}
where
\begin{equation}
\begin{array}[t]{cc}
\mathbf{V_{11}}=\frac{1}{12}(\mathbf{A_{11}^{\mathit{2}}}+\mathbf{A_{12}}\mathbf{A_{21}}); & \mathbf{V_{21}}=\frac{1}{12}(\mathbf{A_{21}}\mathbf{A_{11}}-\mathbf{A_{11}^{\mathrm{T}}}\mathbf{A_{21}})\end{array}
\end{equation}
Note that plate kinematics is then enforced such that $\theta_z=u_x=u_y=0$ to reduce the size of the stiffness matrix. The matrix equation governing free vibrations may be expressed as
\begin{equation}
(\mathbf{K}+\omega^{2}\mathbf{M})\mathbf{v}=\mathbf{0}\label{eq:freevibraiton}
\end{equation}
in which the global mass matrix of the plate structures is defined as
\begin{equation}
\mathbf{M}=\rho\int_{\Omega}\mathbf{N}{}^{\mathrm{T}}\mathbf{H}\mathbf{N}~\rmd \Omega \label{eq:GolbalMass-1}
\end{equation}
with the transformation matrix
\begin{equation}
\renewcommand{\arraystretch}{2}
\mathbf{H}=\left[ \begin{array}{ccc}
\frac{h^{3}}{12} & 0 & 0\\
0 & \frac{h^{3}}{12} & 0\\
0 & 0 & h
\end{array}\right]
\end{equation}
where $\mathbf{H}$ is used to transform the 3D consistent mass matrix obtained from the kinetic energy into the one with plate DOFs, $\omega$ is the natural frequency and $\mathbf{v}$ is the corresponding mode shape. The stability problem involves the solution of the following eigenproblem
\begin{equation}
(\mathbf{K}+\lambda\mathbf{K}_{G})\mathbf{d}=\mathbf{0}\label{eq:buckling}
\end{equation}
where $\lambda$ is the critical buckling load parameter, a constant by which the in-plane loads must be multiplied to cause buckling and $\mathbf{K}_{G}$ is the geometric stiffness matrices. The geometric stiffness matrix is given by:
\begin{equation}
\mathbf{K}_G = \int \left\{ \mathbf{G}_b^{\rm T} \mathbf{N}^\circ \mathbf{G}_b +  \mathbf{G}_{s1}^{\rm T} \mathbf{N}^\circ \mathbf{G}_{s1} + \mathbf{G}_{s2}^{\rm T} \mathbf{N}^\circ \mathbf{G}_{s2} \right\}~\rmd \Omega
\end{equation}
where
\begin{align}
\renewcommand{\arraystretch}{2}
\mathbf{G}_b &= \left[ \begin{array}{ccc} \frac{\partial N}{\partial x} & 0 & 0 \\ \frac{\partial N}{\partial y} & 0 & 0 \end{array} \right]; \nonumber \\
\mathbf{G}_{s1} &= \left[ \begin{array}{ccc} 0 & \frac{\partial N}{\partial x} &  0 \\ 0 & \frac{\partial N}{\partial y} & 0 \end{array} \right]; \nonumber \\
\mathbf{G}_{s2} &= \left[ \begin{array}{ccc} 0 & 0 & \frac{\partial N}{\partial x} \\ 0 & 0 & \frac{\partial N}{\partial y} \end{array} \right].
\end{align}
and
\begin{equation}
\mathbf{N}^\circ = \left[ \begin{array}{cc} N_x & N_{xy} \\ N_{xy} & N_y \end{array} \right]
\end{equation}
The natural frequency and critical buckling load in \Erefs{eq:freevibraiton} and (\ref{eq:buckling}) are computed using a standard eigenvalue algorithm.

\section{Element description}
\label{elemdesc}
\begin{figure}[htpb]
\centering
\subfigure[]{\includegraphics[width=8cm]{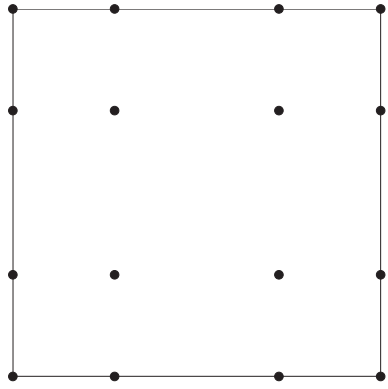}}
\subfigure[]{\includegraphics[width=8cm]{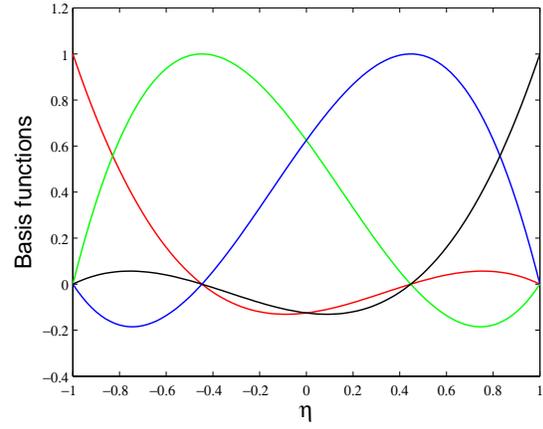}}
\subfigure[]{\includegraphics[width=10cm]{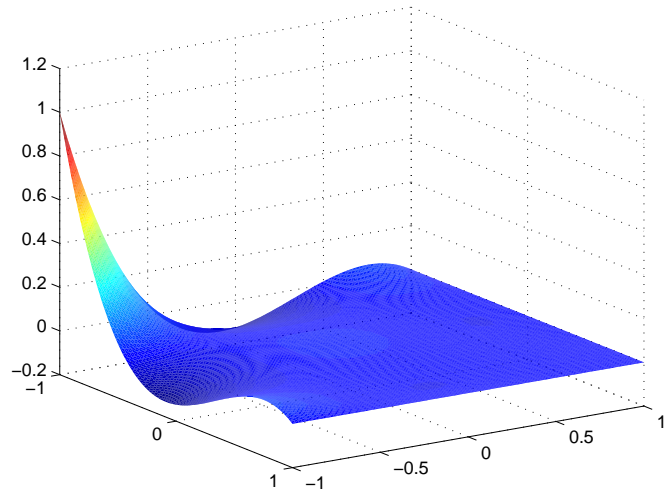}}
\caption{3rd order element: (a) Nodal location (b) 1D  shape functions (c) 2D shape function.}
\label{fig:GLL1D}
\end{figure}
High-order spectral elements are used in this study to discretize the in-plane dimensions of the plate. As the formulation presented here is 3D consistent, in order to represent the constant curvature, the minimum requirement is that the second derivative of the shape function is a constant. Hence, we require a second order element. However, for better convergence and accuracy, in this study we have employed a 3$^{\rm rd}$ order element (see \fref{fig:GLL1D}). The
2D shape functions are obtained by the product of two sets of 1D shape
functions defined separately in local coordinates $\eta$ and $\zeta$:
\begin{equation}
N_{i}(\eta,\,\zeta)=N_{i_{\eta}}(\eta)N_{i_{\zeta}}(\zeta)\label{eq:2Dshapefunc}
\end{equation}
Denoting the orders of the two 1D elements as $p_{\eta}$ and $p_{\zeta}$, respectively, the total number of nodes of the 2D element is equal to $n_{d}=(p_{\eta}+1)(p_{\zeta}+1)$. The local nodal number $i=1,\,2\,,...,n_{d}$ is defined by the nodal numbers $i_{\eta},\, i_{\zeta}$ of the two 1D elements as
\begin{equation}
i=(i_{\zeta}-1)\times(p_{\eta}+1)+i_{\eta}\label{eq:2Dshapefunc-index}
\end{equation}
and the nodal number ascends firstly along $\eta$ direction and then $\zeta$ direction. This also applies to the weights such that
\begin{equation}
w_{i}=w_{i_{\eta}}w_{i_{\zeta}}.\label{eq:2Dweight}
\end{equation}
The coefficient matrices $\mathbf{E}_0, \mathbf{E}_1$ and $\mathbf{E}_2$, the mass matrix $\mathbf{M}$ and geometry stiffness matrix $\mathbf{K}_G$ are computed using Gauss-Lobatto-Legendre quadrature. As demonstrated in \cite{gravenkampman2013}, $\mathbf{E}_0$ becomes a lumped matrix. Consequently, the inversion of $\mathbf{E}_0$ , which is required to compute the stiffness matrix, becomes trivial in Equation \eqref{eq:Z}. This leads us an efficient stiffness matrix construction. With the same quadrature, a diagonal mass matrix is also obtained and the 3$\times$3 submatrix $\mathbf{M}_{\mathit{i}}$ corresponding to the diagonal block of the $i^{th}$ node is expressed as 
\begin{equation}
\mathbf{M}_{\mathit{i}}= \int\limits_\Omega N_{i}(\eta,\,\zeta)\mathbf{H}\rho N_{i}(\eta,\,\zeta)|J|~\mathrm{d} \Omega =w_{i}\mathbf{H}\rho |J| 
\label{eq:finalm0}
\end{equation}

\section{Numerical Results}
\label{numres}
In this section, we present examples of the free vibration and the mechanical buckling of plates with in-plane material inhomogeneity based on the approach discussed in the previous section. The effect of various parameters, viz., material gradient index $n$, the skewness of the plate, $\psi$, the plate aspect ratio $a/b$, the plate thickness $b/h$ and the boundary conditions on the global response is numerically studied. The FGM plate considered in this study is made up of silicon nitride (Si$_3$N$_4$) and stainless steel (SUS304). The material is considered to be temperature independent. The mass density $(\rho)$, Young's modulus $(E)$ and Poisson's ratio $(\nu)$ are $\rho_c=$ 2370 Kg/m$^3$, $E_c=$ 348.43 GPa, $\nu_c=$ 0.24 for Si$_3$N$_4$ and $\rho_m=$ 8166 Kg/m$^3$, $E_m=$ 201.04 GPa, $\nu_m=$ 0.3262 for SUS304.

\paragraph{Validation} Before proceeding with the detailed numerical study, the formulation developed herein is validated against available numerical results pertaining to isotropic plates for fundamental frequency and critical buckling load. Table \ref{table:convergenceRes} compares the first non-dimensionalized fundamental frequency and critical buckling load factor for a simple supported square plate based on a progressive mesh refinement. It is observed that with decreasing element size, the non-dimensionalized frequency and the non-dimensionalized critical buckling load converges. It can be seen that the results from the present formulation compare very well with the available solutions, a structured mesh of 8$\times$8 with 3$^{\rm rd}$ order element is found to be adequate to model the full plate. This mesh is, therefore, used in the subsequent studies of this section. 

\begin{table}[htpb]
\centering
\renewcommand{\arraystretch}{1.2}
\caption{Convergence of the non-dimensionalized fundamental frequency $(\overline{\omega})$ and the critical buckling load parameter $(\lambda_{cru})$ with plate thickness ratio $b/h=1000$. The plate is simply supported on all the edges.}
\begin{tabular}{lrrr}
\hline
Mesh & $\overline{\omega} = \omega \left( \frac{a}{\pi} \right)^2 \sqrt{\frac{\rho_c}{D_c}}$ && $\lambda_{cru} = \frac{N_{xxcr}^{\circ}b^{2}}{\pi^{2}D_{c}} $ \\
\hline
2$\times$2 & 2.0035 &&  4.0154 \\
4$\times$4 & 1.9997 && 4.0010 \\
8$\times$8 & 1.9995 &&  4.0000 \\
Ref.~\cite{reddy1984} & 1.9974 && 4.0000 \\
Ref.~\cite{aydogdu2009} & 1.9974 && - \\
Ref.~\cite{liewhung19993} & 1.9993 && - \\
\hline
\end{tabular}
\label{table:convergenceRes}
\end{table}

\subsection{Free vibration}
\begin{figure}[htpb]
\begin{centering}
\subfigure[$\overline{\omega} _1$]{\scalebox{0.65}{
%
%
\begin{tikzpicture}

\begin{axis}[%
width=4.52083333333333in,
height=3.46354166666667in,
scale only axis,
xmin=0,
xmax=2.5,
xlabel={$a/b$},
ymin=0,
ymax=6,
ylabel={$\overline{\omega}_1$},
legend style={draw=none,legend cell align=left}
]
\addplot [color=black,solid,line width=1.0pt,mark size=6.0pt,mark=diamond,mark options={solid,draw=black}]
  table[row sep=crcr]{0.5	5\\
0.75	2.7778\\
1	2\\
1.25	1.64\\
1.5	1.4444\\
1.75	1.3265\\
2	1.25\\
};
\addlegendentry{$n=$ 0};

\addplot [color=black,solid,line width=1.0pt,mark size=3.5pt,mark=x,mark options={solid,draw=black}]
  table[row sep=crcr]{0.5	2.9797\\
0.75	1.6524\\
1	1.1865\\
1.25	0.9694\\
1.5	0.8501\\
1.75	0.7767\\
2	0.7278\\
};
\addlegendentry{$n=$ 1};

\addplot [color=black,solid,line width=1.0pt,mark size=2.5pt,mark=square,mark options={solid,draw=black}]
  table[row sep=crcr]{0.5	2.2412\\
0.75	1.252\\
1	0.9021\\
1.25	0.7381\\
1.5	0.6479\\
1.75	0.5929\\
2	0.5567\\
};
\addlegendentry{$n=$ 5};

\addplot [color=black,solid,line width=1.0pt,mark size=3.5pt,mark=o,mark options={solid,draw=black}]
  table[row sep=crcr]{0.5	2.1533\\
0.75	1.2022\\
1	0.8665\\
1.25	0.7096\\
1.5	0.6236\\
1.75	0.5713\\
2	0.5372\\
};
\addlegendentry{$n=$ 10};

\end{axis}
\end{tikzpicture}%
}}\hspace{1cm}
\subfigure[$\overline{\omega} _2$]{\scalebox{0.65}{
%
%
\begin{tikzpicture}

\begin{axis}[%
width=4.52083333333333in,
height=3.46354166666667in,
scale only axis,
xmin=0,
xmax=2.5,
xlabel={$a/b$},
ymin=0,
ymax=9,
ylabel={$\overline{\omega}_2$},
legend style={draw=none,legend cell align=left}
]
\addplot [color=black,solid,line width=1.0pt,mark size=6.0pt,mark=diamond,mark options={solid,draw=black}]
  table[row sep=crcr]{0.5	8\\
0.75	5.7778\\
1	5\\
1.25	3.56\\
1.5	2.7778\\
1.75	2.3061\\
2	2\\
};
\addlegendentry{$n=$ 0};

\addplot [color=black,solid,line width=1.0pt,mark size=3.5pt,mark=x,mark options={solid,draw=black}]
  table[row sep=crcr]{0.5	4.7459\\
0.75	3.4002\\
1	2.9111\\
1.25	2.1562\\
1.5	1.684\\
1.75	1.3995\\
2	1.2149\\
};
\addlegendentry{$n=$ 1};

\addplot [color=black,solid,line width=1.0pt,mark size=2.5pt,mark=square,mark options={solid,draw=black}]
  table[row sep=crcr]{0.5	3.6084\\
0.75	2.5917\\
1	2.2269\\
1.25	1.6441\\
1.5	1.2835\\
1.75	1.0655\\
2	0.9237\\
};
\addlegendentry{$n=$ 5};

\addplot [color=black,solid,line width=1.0pt,mark size=3.5pt,mark=o,mark options={solid,draw=black}]
  table[row sep=crcr]{0.5	3.4658\\
0.75	2.4942\\
1	2.1487\\
1.25	1.5568\\
1.5	1.2158\\
1.75	1.0097\\
2	0.8755\\
};
\addlegendentry{$n=$ 10};

\end{axis}
\end{tikzpicture}%
}}\par\end{centering}
\begin{centering}
\subfigure[$\overline{\omega} _3$]{\scalebox{0.65}{
%
%
\begin{tikzpicture}

\begin{axis}[%
width=4.52083333333333in,
height=3.46354166666667in,
scale only axis,
xmin=0,
xmax=2.5,
xlabel={$a/b$},
ymin=0,
ymax=14,
ylabel={$\overline{\omega}_3$},
legend style={draw=none,legend cell align=left}
]
\addplot [color=black,solid,line width=1.0pt,mark size=6.0pt,mark=diamond,mark options={solid,draw=black}]
  table[row sep=crcr]{0.5	12.999\\
0.75	8.1111\\
1	5\\
1.25	4.64\\
1.5	4.4445\\
1.75	3.9386\\
2	3.2498\\
};
\addlegendentry{$n=$ 0};

\addplot [color=black,solid,line width=1.0pt,mark size=3.5pt,mark=x,mark options={solid,draw=black}]
  table[row sep=crcr]{0.5	7.65\\
0.75	4.9051\\
1	3.0258\\
1.25	2.6674\\
1.5	2.5202\\
1.75	2.4025\\
2	1.9844\\
};
\addlegendentry{$n=$ 1};

\addplot [color=black,solid,line width=1.0pt,mark size=2.5pt,mark=square,mark options={solid,draw=black}]
  table[row sep=crcr]{0.5	5.8308\\
0.75	3.7366\\
1	2.3067\\
1.25	2.0539\\
1.5	1.9574\\
1.75	1.8388\\
2	1.5179\\
};
\addlegendentry{$n=$ 5};

\addplot [color=black,solid,line width=1.0pt,mark size=3.5pt,mark=o,mark options={solid,draw=black}]
  table[row sep=crcr]{0.5	5.6115\\
0.75	3.5341\\
1	2.1831\\
1.25	1.9869\\
1.5	1.8983\\
1.75	1.7364\\
2	1.4333\\
};
\addlegendentry{$n=$ 10};

\end{axis}
\end{tikzpicture}%
}}\hspace{1cm}
\subfigure[$\overline{\omega} _4$]{\scalebox{0.65}{
%
%
\begin{tikzpicture}

\begin{axis}[%
width=4.47708333333333in,
height=3.33854166666667in,
scale only axis,
xmin=0,
xmax=2.5,
xlabel={$a/b$},
ymin=0,
ymax=18,
ylabel={$\overline{\omega} _4$},
legend style={draw=none,legend cell align=left}
]
\addplot [color=black,solid,line width=1.0pt,mark size=6.0pt,mark=diamond,mark options={solid,draw=black}]
  table[row sep=crcr]{0.5	16.9999\\
0.75	10.7774\\
1	8\\
1.25	6.56\\
1.5	4.9998\\
1.75	4.3266\\
2	4.25\\
};
\addlegendentry{$n=$ 0};

\addplot [color=black,solid,line width=1.0pt,mark size=3.5pt,mark=x,mark options={solid,draw=black}]
  table[row sep=crcr]{0.5	10.2751\\
0.75	6.2362\\
1	4.8597\\
1.25	3.9912\\
1.5	3.047\\
1.75	2.4193\\
2	2.3445\\
};
\addlegendentry{$n=$ 1};

\addplot [color=black,solid,line width=1.0pt,mark size=2.5pt,mark=square,mark options={solid,draw=black}]
  table[row sep=crcr]{0.5	7.8194\\
0.75	4.7846\\
1	3.6948\\
1.25	3.0241\\
1.5	2.3331\\
1.75	1.8972\\
2	1.8568\\
};
\addlegendentry{$n=$ 5};

\addplot [color=black,solid,line width=1.0pt,mark size=3.5pt,mark=o,mark options={solid,draw=black}]
  table[row sep=crcr]{0.5	7.391\\
0.75	4.6226\\
1	3.502\\
1.25	2.8677\\
1.5	2.2029\\
1.75	1.8444\\
2	1.809\\
};
\addlegendentry{$n=$ 10};

\end{axis}
\end{tikzpicture}%
}}\par\end{centering}
\caption{Influence of the plate aspect ratio on the first four non-dimensionalized frequency for a simply supported FGM plate with material gradient indexes $n=0,1,5,10$.}
\label{fig:FreeVibrationSS}
\end{figure}
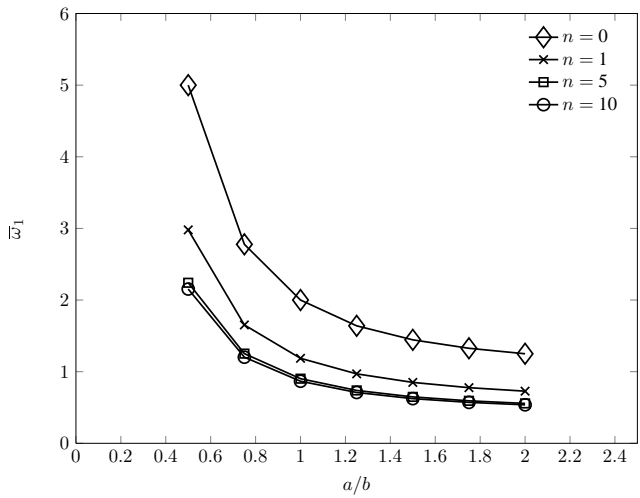
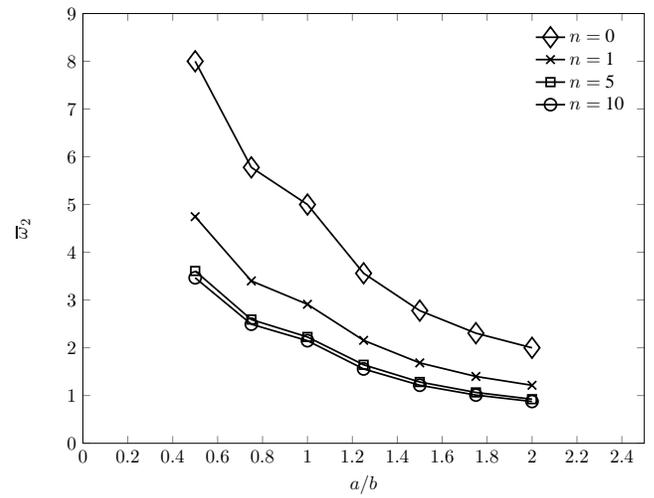
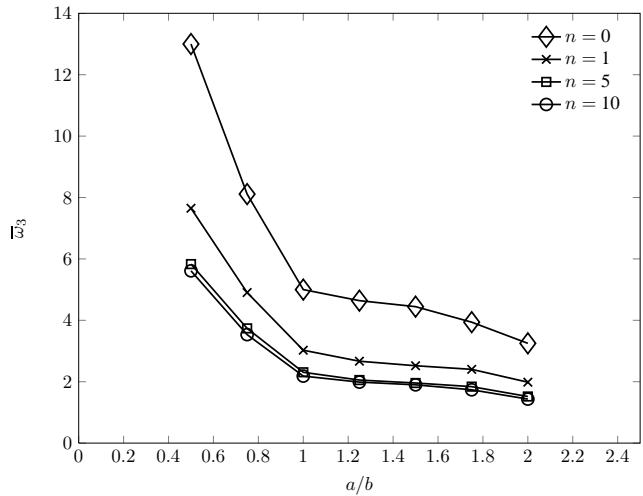
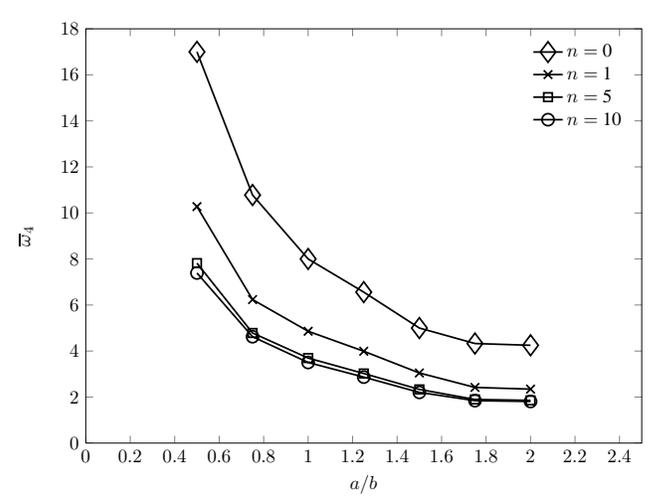

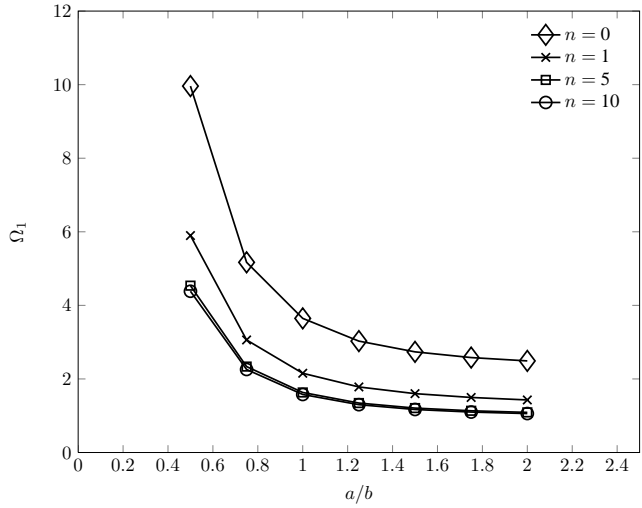
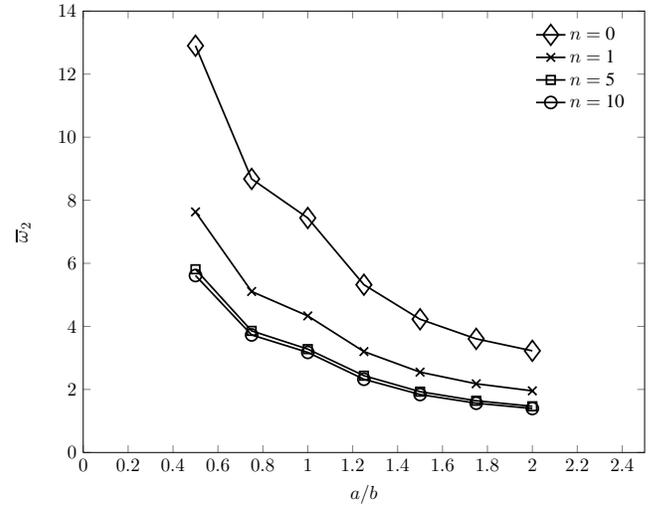
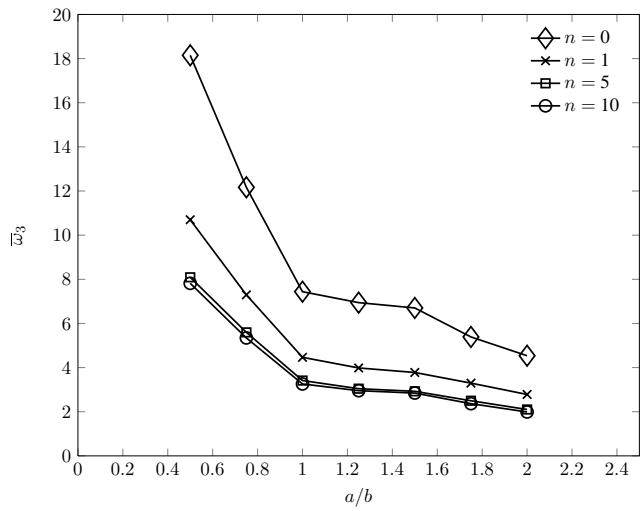
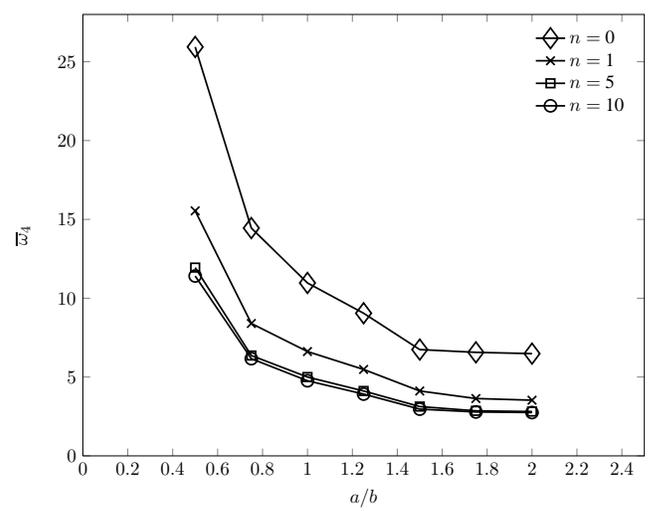
\begin{figure}[htpb]
\begin{centering}
\subfigure[$\overline{\omega} _1$]{\scalebox{0.65}{
%
%
\begin{tikzpicture}

\begin{axis}[%
width=4.52083333333333in,
height=3.55729166666667in,
scale only axis,
xmin=0,
xmax=2.5,
xlabel={$a/b$},
ymin=0,
ymax=12,
ylabel={$\Omega_1$},
legend style={draw=none,legend cell align=left}
]
\addplot [color=black,solid,line width=1.0pt,mark size=6.0pt,mark=diamond,mark options={solid,draw=black}]
  table[row sep=crcr]{0.5	9.961\\
0.75	5.1698\\
1	3.6461\\
1.25	3.0283\\
1.5	2.7362\\
1.75	2.5809\\
2	2.4903\\
};
\addlegendentry{$n=$ 0};

\addplot [color=black,solid,line width=1.0pt,mark size=3.5pt,mark=x,mark options={solid,draw=black}]
  table[row sep=crcr]{0.5	5.8974\\
0.75	3.0597\\
1	2.1543\\
1.25	1.7823\\
1.5	1.6002\\
1.75	1.4964\\
2	1.4288\\
};
\addlegendentry{$n=$ 1};

\addplot [color=black,solid,line width=1.0pt,mark size=2.5pt,mark=square,mark options={solid,draw=black}]
  table[row sep=crcr]{0.5	4.5354\\
0.75	2.3359\\
1	1.632\\
1.25	1.3446\\
1.5	1.2077\\
1.75	1.1342\\
2	1.0906\\
};
\addlegendentry{$n=$ 5};

\addplot [color=black,solid,line width=1.0pt,mark size=3.5pt,mark=o,mark options={solid,draw=black}]
  table[row sep=crcr]{0.5	4.3834\\
0.75	2.257\\
1	1.5773\\
1.25	1.3005\\
1.5	1.1694\\
1.75	1.0996\\
2	1.0589\\
};
\addlegendentry{$n=$ 10};

\end{axis}
\end{tikzpicture}%
}}\hspace{1cm}
\subfigure[$\overline{\omega} _2$]{\scalebox{0.65}{
%
%
\begin{tikzpicture}

\begin{axis}[%
width=4.52083333333333in,
height=3.55729166666667in,
scale only axis,
xmin=0,
xmax=2.5,
xlabel={$a/b$},
ymin=0,
ymax=14,
ylabel={$\overline{\omega}_2$},
legend style={draw=none,legend cell align=left}
]
\addplot [color=black,solid,line width=1.0pt,mark size=6.0pt,mark=diamond,mark options={solid,draw=black}]
  table[row sep=crcr]{0.5	12.901\\
0.75	8.6765\\
1	7.4381\\
1.25	5.3215\\
1.5	4.2264\\
1.75	3.6041\\
2	3.2253\\
};
\addlegendentry{$n=$ 0};

\addplot [color=black,solid,line width=1.0pt,mark size=3.5pt,mark=x,mark options={solid,draw=black}]
  table[row sep=crcr]{0.5	7.6317\\
0.75	5.107\\
1	4.3321\\
1.25	3.2012\\
1.5	2.5492\\
1.75	2.1801\\
2	1.956\\
};
\addlegendentry{$n=$ 1};

\addplot [color=black,solid,line width=1.0pt,mark size=2.5pt,mark=square,mark options={solid,draw=black}]
  table[row sep=crcr]{0.5	5.8115\\
0.75	3.8576\\
1	3.2789\\
1.25	2.4357\\
1.5	1.9295\\
1.75	1.6415\\
2	1.4658\\
};
\addlegendentry{$n=$ 5};

\addplot [color=black,solid,line width=1.0pt,mark size=3.5pt,mark=o,mark options={solid,draw=black}]
  table[row sep=crcr]{0.5	5.613\\
0.75	3.7289\\
1	3.1752\\
1.25	2.3234\\
1.5	1.8393\\
1.75	1.5639\\
2	1.396\\
};
\addlegendentry{$n=$ 10};

\end{axis}
\end{tikzpicture}%
}}\par\end{centering}
\begin{centering}
\subfigure[$\overline{\omega} _3$]{\scalebox{0.65}{
%
%
\begin{tikzpicture}

\begin{axis}[%
width=4.52083333333333in,
height=3.55729166666667in,
scale only axis,
xmin=0,
xmax=2.5,
xlabel={$a/b$},
ymin=0,
ymax=20,
ylabel={$\overline{\omega}_3$},
legend style={draw=none,legend cell align=left}
]
\addplot [color=black,solid,line width=1.0pt,mark size=6.0pt,mark=diamond,mark options={solid,draw=black}]
  table[row sep=crcr]{0.5	18.1493\\
0.75	12.1685\\
1	7.4381\\
1.25	6.9431\\
1.5	6.7017\\
1.75	5.3883\\
2	4.5374\\
};
\addlegendentry{$n=$ 0};

\addplot [color=black,solid,line width=1.0pt,mark size=3.5pt,mark=x,mark options={solid,draw=black}]
  table[row sep=crcr]{0.5	10.6978\\
0.75	7.2949\\
1	4.4653\\
1.25	3.9832\\
1.5	3.7783\\
1.75	3.2964\\
2	2.783\\
};
\addlegendentry{$n=$ 1};

\addplot [color=black,solid,line width=1.0pt,mark size=2.5pt,mark=square,mark options={solid,draw=black}]
  table[row sep=crcr]{0.5	8.0905\\
0.75	5.5952\\
1	3.413\\
1.25	3.0437\\
1.5	2.9252\\
1.75	2.4978\\
2	2.1016\\
};
\addlegendentry{$n=$ 5};

\addplot [color=black,solid,line width=1.0pt,mark size=3.5pt,mark=o,mark options={solid,draw=black}]
  table[row sep=crcr]{0.5	7.8171\\
0.75	5.3432\\
1	3.2577\\
1.25	2.9539\\
1.5	2.846\\
1.75	2.3654\\
2	1.9885\\
};
\addlegendentry{$n=$ 10};

\end{axis}
\end{tikzpicture}%
}}\hspace{1cm}
\subfigure[$\overline{\omega} _4$]{\scalebox{0.65}{
%
%
\begin{tikzpicture}

\begin{axis}[%
width=4.52083333333333in,
height=3.55729166666667in,
scale only axis,
xmin=0,
xmax=2.5,
xlabel={$a/b$},
ymin=0,
ymax=28,
ylabel={$\overline{\omega}_4$},
legend style={draw=none,legend cell align=left}
]
\addplot [color=black,solid,line width=1.0pt,mark size=6.0pt,mark=diamond,mark options={solid,draw=black}]
  table[row sep=crcr]{0.5	25.9384\\
0.75	14.4493\\
1	10.9685\\
1.25	9.0471\\
1.5	6.7424\\
1.75	6.5672\\
2	6.4849\\
};
\addlegendentry{$n=$ 0};

\addplot [color=black,solid,line width=1.0pt,mark size=3.5pt,mark=x,mark options={solid,draw=black}]
  table[row sep=crcr]{0.5	15.5416\\
0.75	8.3983\\
1	6.6187\\
1.25	5.4827\\
1.5	4.116\\
1.75	3.6375\\
2	3.5323\\
};
\addlegendentry{$n=$ 1};

\addplot [color=black,solid,line width=1.0pt,mark size=2.5pt,mark=square,mark options={solid,draw=black}]
  table[row sep=crcr]{0.5	11.9439\\
0.75	6.3655\\
1	5.0008\\
1.25	4.1163\\
1.5	3.1286\\
1.75	2.8559\\
2	2.8108\\
};
\addlegendentry{$n=$ 5};

\addplot [color=black,solid,line width=1.0pt,mark size=3.5pt,mark=o,mark options={solid,draw=black}]
  table[row sep=crcr]{0.5	11.4092\\
0.75	6.165\\
1	4.7654\\
1.25	3.9204\\
1.5	2.965\\
1.75	2.7856\\
2	2.7484\\
};
\addlegendentry{$n=$ 10};

\end{axis}
\end{tikzpicture}%
}}\par\end{centering}
\caption{Influence of the plate aspect ratio on the first four non-dimensionalized frequency for a clamped FGM plate with material gradient indexes $n=0,1,5,10$.} 
\label{fig:FreeVibrationCC}
\end{figure}

In this section, the free vibration of FGM rectangular and skew plate with side lengths $a$ and $b$ and thickness $h$ is studied. In all cases, we present the non-dimenisonalized flexural frequencies, unless specified otherwise, as:
\begin{equation}
\overline{\omega} = \omega \left( \frac{a}{\pi} \right)^2 \sqrt{\frac{\rho_c}{D_c}}
\end{equation}
where $D_c=\frac{E_c h^3}{12(1-\nu^2)}$ and $\rho_c$ are the flexural rigidity and the mass density of the ceramic phase, respectively. The influence of the plate thickness ratio $b/h$, the material gradient index $n$ and the boundary conditions on the first non-dimensionalized frequency is shown in Table \ref{table:effectofahVib}. It indicates that with increasing plate thickness ratio, the non-dimensionalized frequency decreases irrespective of the boundary conditions of the plate. This can be attributed to the decreased flexural rigidity. Table \ref{table:effectofahVib} also shows that the non-dimensionalized frequency reduces with increasing material gradient index $n$ due to the increase in metallic volume fraction. \frefs{fig:FreeVibrationSS} - (\ref{fig:FreeVibrationCC}) illustrates the influence of the plate aspect ratio $a/b$ and the material gradient index $n$ on the first four non-dimensionalized frequency for a plate with simply supported edges and clamped edges, respectively. The non-dimensionalized frequency decreases for both simply supported and clamped boundary conditions with increasing plate aspect ratio $a/b$ ratio or increasing material gradient index $n$. Table \ref{table:effectskewVib} shows the influence of the skew angle $\psi$ and the gradient index $n$ on the first four non-dimensionalized frequency for a square FGM plate with $b/h=$ 100. The plate is subjected to simply support and clamped boundary conditions on alternating edge, i.e, SCSC. With increasing skew angle $\psi$, the non-dimensionalized frequency increases, whilst the frequency decreases with increasing material gradient index $n$. This can be attributed to the change in the flexural rigidity due to the change in the geometry of the plate and the increase in the metallic volume fraction.

\begin{table}[htpb]
\renewcommand{\arraystretch}{1.5}
\centering
\caption{Influence of the boundary conditions, the plate thickness ratio $b/h$ and the material gradient index $n$ on the non-dimensionalized fundamental frequency of a square plate.}
\begin{tabular}{crrrrrrr}
\hline 
 & $b/h$ & \multicolumn{6}{c}{Gradient index, $n$}\\
 \cline{3-8}
 & & 0 & 0.5 & 1 & 2 & 5 & 10 \\
 \hline
\multirow{3}{*}{SSSS} & 1000  & 2.0000  & 1.3976 & 1.1865 & 1.0212 & 0.9021 & 0.8665 \\
 & 100  & 1.9995  & 1.3972 & 1.1861 & 1.0209 & 0.9018 & 0.8662 \\
 & 50  & 1.9978  & 1.3960 & 1.1851 & 1.0200 & 0.9011 & 0.8655 \\
\cline{2-8}
\multirow{3}{*}{CCCC} & 1000  & 3.6461 & 2.5437 & 2.1543  & 1.8470  & 1.6320  & 1.5773 \\
 & 100  & 3.6431 & 2.5415  & 2.1524  & 1.8454  & 1.6306  & 1.5759 \\
 & 50  & 3.6341 & 2.5349  & 2.1467  & 1.8406  & 1.6263  & 1.5718 \\
\hline 
\end{tabular}
\label{table:effectofahVib}
\end{table}

\begin{table}[htpb]
\centering
\renewcommand{\arraystretch}{1.2}
\caption{Influence of the skew angle $\psi$ and the material gradient index $n$ on the first four non-dimensionalized frequencies of a SCSC square FGM plate with $b/h=$ 100.}
\begin{tabular}{ccrrrrrrr}
\hline 
 & Skew & \multicolumn{7}{c}{Gradient index, $n$}  \\
\cline{3-9} 
 & Angle & \multicolumn{2}{c}{0} & 0.5 & 1 & 2 & 5 & 10\\
\cline{3-4} 
 & $\psi$ & Ref.~\cite{Lee2004} & Present  \\
\hline 
\multirow{4}{*}{$\overline{\omega} _1$} & $0^{\circ}$ & 2.934 & 2.9316 & 2.0419 & 1.7352 & 1.4960 & 1.3275 & 1.2818\\
 & $15^\circ$ & 3.111 & 3.1107 & 2.3818 & 2.0455 & 1.7445 & 1.5076 & 1.4516\\
 & $30^\circ$ & 3.746 & 3.7555 & 3.1157 & 2.7422 & 2.3501 & 1.9916 & 1.9151\\
 & $45^\circ$ & 5.341 & 5.3647 & 4.7620 & 4.3136 & 3.7529 & 3.1333 & 2.9658\\
\hline 
\multirow{4}{*}{$\overline{\omega} _2$} & $0^\circ$ & 5.548 & 5.5466 & 3.8578 & 3.2565 & 2.7880 & 2.4626 & 2.3795\\
 & $15^\circ$ & 5.765 & 5.7502 & 4.3782 & 3.7537 & 3.1938 & 2.7414 & 2.6146\\
 & $30^\circ$ & 6.514 & 6.5364 & 5.3893 & 4.7608 & 4.1176 & 3.5244 & 3.3750\\
 & $45^\circ$ & 8.488 & 8.5580 & 7.5554 & 6.9019 & 6.1311 & 3.1263 & 2.9457\\
\hline 
\multirow{4}{*}{$\overline{\omega} _3$} & $0^\circ$ & 7.024 & 7.0242 & 4.9013 & 4.2054 & 3.6514 & 3.2228 & 3.0785\\
 & $15^\circ$ & 7.579 & 7.5625 & 5.8082 & 5.0302 & 4.3229 & 3.7304 & 3.5363\\
 & $30^\circ$ & 8.450 & 9.4358 & 7.8143 & 6.8921 & 5.9187 & 4.9814 & 4.6545\\
 & $45^\circ$ & 12.559 & 12.6199 & 11.0078 & 10.0606 & 8.9901 & 5.3264 & 5.1068\\
\hline 
\multirow{4}{*}{$\overline{\omega} _4$} & $0^\circ$ & 9.586 & 9.5671 & 6.7129 & 5.7526 & 4.9765 & 4.3725  & 4.1708\\
 & $15^\circ$ & 9.552 & 9.5381 & 7.2636 & 6.2959 & 5.4370 & 4.7372 & 4.5186\\
 & $30^\circ$ & 10.224 & 10.2492 & 8.4668 & 7.5816 & 6.6930 & 5.9038 & 5.7063\\
 & $45^\circ$ & 13.899 & 13.9775 & 12.5132 & 11.4675 & 10.1024 & 7.6899 & 7.0433\\
\hline 
\end{tabular}
\label{table:effectskewVib}
\end{table}

\subsection{Buckling}
\begin{figure}[hp]
\begin{centering}
\subfigure[]{\scalebox{0.7}{
%
%
\begin{tikzpicture}

\begin{axis}[%
width=4.52083333333333in,
height=3.55729166666667in,
scale only axis,
xmin=0,
xmax=2.5,
xlabel={$a/b$},
ymin=0,
ymax=8,
ytick={0, 2, 4, 6, 8},
ylabel={$\lambda_{cru}$},
legend style={draw=none,legend cell align=left}
]
\addplot [color=black,solid,line width=1.0pt,mark size=6.0pt,mark=diamond,mark options={solid,draw=black}]
  table[row sep=crcr]{0.5	6.2501\\
0.75	4.3403\\
1	4\\
1.25	4.2026\\
1.5	4.3412\\
1.75	4.0724\\
2	4.0005\\
};
\addlegendentry{$n=$ 0};

\addplot [color=black,solid,line width=1.0pt,mark size=3.5pt,mark=x,mark options={solid,draw=black}]
  table[row sep=crcr]{0.5	4.9909\\
0.75	3.4577\\
1	3.1651\\
1.25	3.2439\\
1.5	3.2305\\
1.75	3.0843\\
2	3.0069\\
};
\addlegendentry{$n=$ 1};

\addplot [color=black,solid,line width=1.0pt,mark size=2.5pt,mark=square,mark options={solid,draw=black}]
  table[row sep=crcr]{0.5	4.093\\
0.75	2.8721\\
1	2.6435\\
1.25	2.73\\
1.5	2.7627\\
1.75	2.6378\\
2	2.583\\
};
\addlegendentry{$n=$ 5};

\addplot [color=black,solid,line width=1.0pt,mark size=3.5pt,mark=o,mark options={solid,draw=black}]
  table[row sep=crcr]{0.5	3.9404\\
0.75	2.7623\\
1	2.5474\\
1.25	2.6541\\
1.5	2.7157\\
1.75	2.5728\\
2	2.5229\\
};
\addlegendentry{$n=$10};

\end{axis}
\end{tikzpicture}%
}}\hspace{1cm}
\subfigure[]{\scalebox{0.7}{
%
%
\begin{tikzpicture}

\begin{axis}[%
width=4.52083333333333in,
height=3.55729166666667in,
scale only axis,
xmin=0,
xmax=2.5,
xlabel={$a/b$},
ymin=0,
ymax=6,
ytick={0, 2, 4, 6},
ylabel={$\lambda_{crb}$},
legend style={draw=none,legend cell align=left}
]
\addplot [color=black,solid,line width=1.0pt,mark size=6.0pt,mark=diamond,mark options={solid,draw=black}]
  table[row sep=crcr]{0.5	5\\
0.75	2.7778\\
1	2\\
1.25	1.64\\
1.5	1.4445\\
1.75	1.3266\\
2	1.25\\
};
\addlegendentry{$n=$ 0};

\addplot [color=black,solid,line width=1.0pt,mark size=3.5pt,mark=x,mark options={solid,draw=black}]
  table[row sep=crcr]{0.5	3.9962\\
0.75	2.2192\\
1	1.5962\\
1.25	1.3069\\
1.5	1.1488\\
1.75	1.0526\\
2	0.9893\\
};
\addlegendentry{$n=$ 1};

\addplot [color=black,solid,line width=1.0pt,mark size=2.5pt,mark=square,mark options={solid,draw=black}]
  table[row sep=crcr]{0.5	3.2758\\
0.75	1.8411\\
1	1.3285\\
1.25	1.0856\\
1.5	0.9507\\
1.75	0.8678\\
2	0.8131\\
};
\addlegendentry{$n=$ 5};

\addplot [color=black,solid,line width=1.0pt,mark size=3.5pt,mark=o,mark options={solid,draw=black}]
  table[row sep=crcr]{0.5	3.1526\\
0.75	1.7688\\
1	1.2761\\
1.25	1.0437\\
1.5	0.9153\\
1.75	0.8367\\
2	0.7852\\
};
\addlegendentry{$n=$ 10};

\end{axis}
\end{tikzpicture}%
}}\par\end{centering}
\caption{Critical buckling parameters of a rectangular plate with material gradient indexes $n=0,1,5,10$ : (a) uni-axial compressive load; (b) bi-axial compressive load.}
\label{fig:bucklingRectangular}
\end{figure}
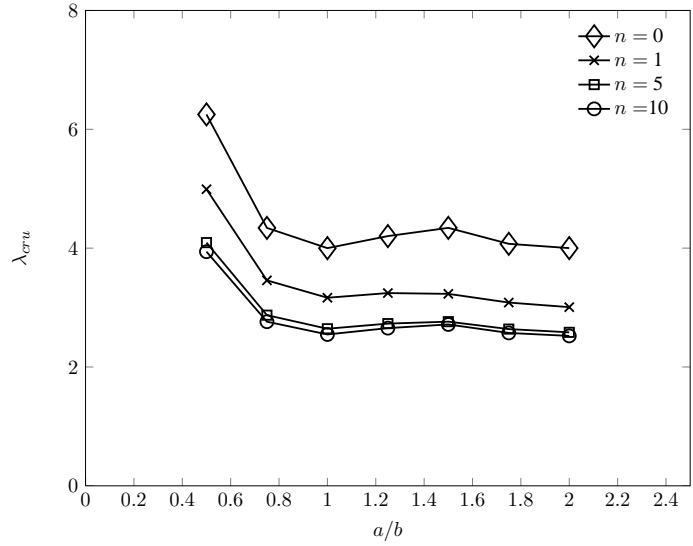
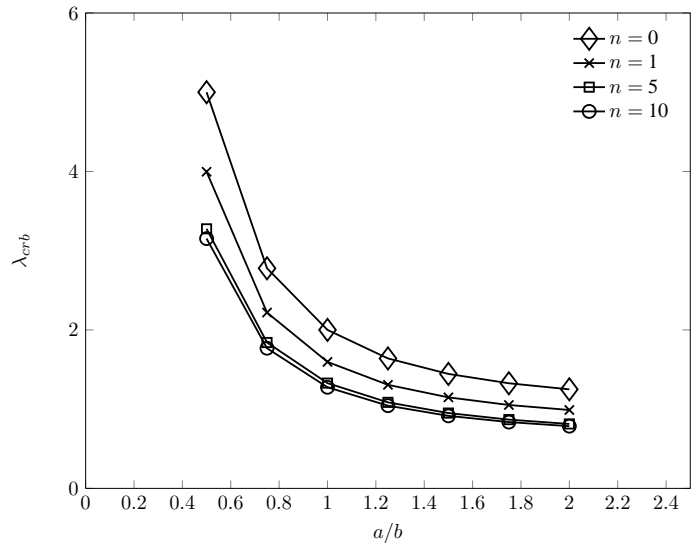

In this section, we present the mechanical buckling behavior of FGM rectangular and skew plates with in-plane material inhomogeneity under uni- and bi-axial compressive loads. In all our cases, we present the non-dimensionalized critical buckling parameters, unless otherwise specified, as:
\begin{align}
\lambda_{cru} & = \frac{N_{xxcr}^{\circ}b^{2}}{\pi^{2}D_{c}} \nonumber \\
\lambda_{crb} & = \frac{N_{yycr}^{\circ}b^{2}}{\pi^{2}D_{c}}
\end{align}
where $\lambda_{cru}$ and $\lambda_{crb}$ are the critical buckling parameters corresponding to uni- and bi-axial compressive loads, $D_{c}=\frac{E_{c}h^{3}}{12(1-v^{2})}$. The influence of plate aspect ratio $a/b$ on the critical buckling load parameter is shown in \fref{fig:bucklingRectangular} subjected to uni- and bi-axial compressive loads. The influence of material gradient index $n$ is also shown. It is seen that the critical buckling parameter decreases with increasing plate aspect ratio and material gradient index. 

The influence of the plate thickness ratio $b/h$, the material gradient index $n$ and boundary conditions on the critical buckling parameter is given in Table \ref{table:effectahBuck}. The influence of the skew angle $\psi$ on the critical buckling parameters of a simply supported FGM skew plate with $b/h=$ 1000 is given in Table \ref{table:effectskewBuck}. It can be seen that the critical buckling parameters obtained from proposed technique are in good agreement with the results available in the literature~\cite{ganapathiprakash2006} for an isotropic plate. It is seen that increasing the critical buckling parameter for both uni- and bi-axial compressive load decreases with increasing material gradient index and the skew angle. This can be attributed to the decreasing flexural rigidity and increasing metallic volume fraction respectively.

\begin{table}[htpb]
\renewcommand{\arraystretch}{1.5}
\centering
\caption{Influence of the plate thickness ratio $b/h$, the material gradient index $n$ and the boundary conditions on the critical buckling parameters for a square FGM plate.}
\begin{tabular}{crrrrrrrrr}
\hline
 & $b/h$ &  & \multicolumn{7}{c}{Gradient index, $n$}\\
\cline{4-10}
 &  &  & \multicolumn{2}{c}{0} & 0.5 & 1 & 2 & 5 & 10 \\
\cline{4-5}
 &  &  & Ref.~\cite{wang1997} & Present \\
\hline
\multirow{6}{*}{SSSS} & \multirow{2}{*}{1000} & $\lambda_{cru}$ & 4.0000 & 4.0000 & 3.4553 & 3.1651  & 2.8861 & 2.6435 & 2.5474\\
 &  & $\lambda_{crb}$ & 2.0000 & 2.0000 & 1.7356 & 1.5962  & 1.4573 & 1.3285 & 1.2761\\
 & \multirow{2}{*}{100} & $\lambda_{cru}$ &  & 3.9979 & 3.4533 & 3.1632 & 2.8844 & 2.6419 & 2.5459\\
 &  & $\lambda_{crb}$ &  & 1.9989 & 1.7346 & 1.5953 & 1.4565 & 1.3277 & 1.2754\\
 & \multirow{2}{*}{50} & $\lambda_{cru}$ &  & 3.9913 & 3.4473 & 3.1574 & 2.8791 & 2.6372 & 2.5414\\
 &  & $\lambda_{crb}$ &  & 1.9957 & 1.7316 & 1.5925 & 1.4540 & 1.3254 & 1.2732\\
\hline 
\multirow{6}{*}{CCCC} & \multirow{2}{*}{1000} & $\lambda_{cru}$ & 10.0740 & 10.0782 & 8.7605 & 7.9632 & 7.1523 & 6.4679 & 6.2671\\
 &  & $\lambda_{crb}$ & 5.3036 & 5.3054 & 4.6056 & 4.2238 & 3.8288 & 3.4712 & 3.3533\\
 & \multirow{2}{*}{100} & $\lambda_{cru}$ &  & 10.0624 & 8.7426 & 7.9463 & 7.1369 & 6.4541 & 6.2538\\
 &  & $\lambda_{crb}$ &  & 5.2968 & 4.5980 & 4.2167 & 3.8222 & 3.4652 & 3.3475\\
 & \multirow{2}{*}{50} & $\lambda_{cru}$ &  & 10.0034 & 8.6894 & 7.8962 & 7.0910 & 6.4129 & 6.2143\\
 &  & $\lambda_{crb}$ &  & 5.2711 & 4.5752 & 4.1954 & 3.8026 & 3.4473 & 3.3301\\
\hline 
\end{tabular}
\label{table:effectahBuck}
\end{table}

\begin{table}[htpb]
\renewcommand{\arraystretch}{1.2}
\centering
\caption{Influence of the skew angle $\psi$ and the material gradient index $n$ on the critical buckling parameters for a simply supported in-plan FGM skew square plate with $b/h=$ 1000.}
\begin{tabular}{ccccccccc}
\hline 
Skew &  & \multicolumn{7}{c}{Gradient index, $n$} \\
\cline{3-9} 
Angle& & \multicolumn{2}{c}{0} & 0.5 & 1 & 2 & 5 & 10 \\
\cline{3-4}
$\psi$ & & Ref.~\cite{ganapathiprakash2006} & Present \\
\hline 
\multirow{2}{*}{$0^{\circ}$} & $\lambda_{cru}$ & 4.0000 & 4.0000 & 3.4553 & 3.1651  & 2.8861 & 2.6435 & 2.5474\tabularnewline
 & $\lambda_{crb}$ & 2.0000 & 2.0000 & 1.7356 & 1.5962  & 1.4573 & 1.3285 & 1.2761\tabularnewline
\multirow{2}{*}{$15^{\circ}$} & $\lambda_{cru}$ & 4.3946 & 4.4026 & 3.9951 & 3.7161  & 3.3932 & 3.0937 & 3.0435\tabularnewline
 & $\lambda_{crb}$ & 2.1154 & 2.1191 & 1.9267 & 1.7997  & 1.6537 & 1.5143 & 1.4991\tabularnewline
\multirow{2}{*}{$30^{\circ}$} & $\lambda_{cru}$ & 5.8966 & 5.9316 & 5.5845 & 5.3066  & 4.9138 & 4.4514 & 4.4811\tabularnewline
 & $\lambda_{crb}$ & 2.5365 & 2.5503 & 2.4006 & 2.2868  & 2.1344 & 1.9766 & 2.0898\tabularnewline
\multirow{2}{*}{$45^{\circ}$} & $\lambda_{cru}$ & 10.1031 & 10.1171 & 9.7879 & 9.4940  & 9.0046 & 8.1707 & 8.2645\tabularnewline
 & $\lambda_{crb}$ & 3.6399 & 3.6326 & 3.5109 & 3.4063  & 3.2425 & 3.0212 & 3.5359\tabularnewline
\hline 
\end{tabular}
\label{table:effectskewBuck}
\end{table}

\section{Conclusions}
In this article, we presented a three dimensional consistent approach that does not require ad hoc shear correction factors to analyse plate structures. Based on this approach, we studied the free vibration and mechanical buckling of thin functionally graded material plates considering various parameters such as the material gradient index, the thickness ratio, the plate aspect ratio and the boundary conditions. From the detailed numerical study, it can be concluded that the material gradient index has strong influence on the fundamental frequency and the critical buckling load parameter. It is also observed that the change in the non-dimensionalized frequency and the critical buckling load parameter is significant for material gradient index $n \le 2$.

\section*{Acknowledgements} 
Sundararajan Natarajan would like to acknowledge the financial support of the School of Civil and Environmental Engineering, The University of New South Wales for his research fellowship since September 2012. 

\section*{References}
\bibliographystyle{plain}
\bibliography{fgminplane} 

\begin{thebibliography}{10}

\bibitem{aydogdu2009}
M~Aydogdu.
\newblock A new shear deformation theory for laminated composite plates.
\newblock {\em Composite Structures}, 89:94--101, 2009.

\bibitem{bathedvorkin1985}
KJ~Bathe and EN~Dvorkin.
\newblock A four node plate bending element based on {Mindlin/Reissner} plate
  theory and a mixed interpolation.
\newblock {\em International Journal for Numerical Methods in Engineering},
  21:367--383, 1985.

\bibitem{brezzibathe1989}
F~Brezzi, KJ~Bathe, and M~Fortin.
\newblock Mixed interpolated elements for {Reissner-Mindlin} plates.
\newblock {\em International Journal for Numerical Methods in Engineering},
  28:1787--1801, 1989.

\bibitem{carrera2003}
E~Carrera.
\newblock Theories and finite elements for multilayered plates and shells: A
  unified compact formulation with numerical assessment and benchmarking.
\newblock {\em Archives of Computational Methods in Engineering}, 10:215--296,
  2003.

\bibitem{carrerademasi2002}
E~Carrera and L~Demasi.
\newblock Classical and advanced multilayered plate elements based upon {PVD}
  and {RMVT. Part 1: derivation of finite element matrices}.
\newblock {\em International Journal for Numerical Methods in Engineering},
  55:191--231, 2002.

\bibitem{ferreirabatra2005}
AJM Ferreira, RC~Batra, CMC Roque, LF~Qian, and PALS Martins.
\newblock Static analysis of functionally graded plates using third-order shear
  deformation theory and a meshless method.
\newblock {\em Composite Structures}, 69:449--457, 2005.

\bibitem{ferreirabatra2006}
AJM Ferreira, RC~Batra, CMC Roque, LK~Qian, and RMN Jorge.
\newblock Natural frequencies of functionally graded plates by a meshless
  method.
\newblock {\em Composite Structures}, 75:593--600, 2006.

\bibitem{ganapathiprakash2006}
M~Ganapathi, T~Prakash, and N~Sundararajan.
\newblock Influence of functionally graded material on buckling of skew plates
  under mechanical loads.
\newblock {\em Journal of Engineering Mechanics - {ASCE}}, 132:902--905, 2006.

\bibitem{goupeevel2006}
AJ~Goupee and SS~Vel.
\newblock Optimization of natural frequencies of bidirectional functionally
  graded beams.
\newblock {\em Struct Multidiscip Optim}, 32:473--484, 2006.

\bibitem{gravenkampman2013}
H~Gravenkamp, H~Man, C~Song, and J~Prager.
\newblock The computation of dispersion relations for three-dimensional elastic
  waveguides using the scaled boundary finite element method.
\newblock {\em Journal of Sound and Vibration}, 332:3756--3771, 2013.

\bibitem{greimannlynn1970}
LF~Greimann and PP~Lynn.
\newblock Finite element analysis of plate bending with transverse shear
  deformation.
\newblock {\em Nuclear Engineering and Design}, 14:223--230, 1970.

\bibitem{hughescohen1978}
TJF Hughes, M~Cohen, and M~Haroun.
\newblock Reduced and selective integration technique in finite element method
  of plates.
\newblock {\em Nuclear Engineering and Design}, 46:203--222, 1978.

\bibitem{jhakant2013}
DK~Jha, Tarun Kant, and RK~Singh.
\newblock A critical review of recent research on functionally graded plates.
\newblock {\em Composite Structures}, 96:833--849, 2013.

\bibitem{Lee2004}
S.J. Lee.
\newblock {Free vibration analysis of plates by using a four-node finite
  element formulated with assumed natural transverse shear strain}.
\newblock {\em Journal of Sound and Vibration}, 278:657--684, 2004.

\bibitem{liewhung19993}
KM~Liew, KC~Hung, and KM~Lim.
\newblock A continuum three-dimensional vibration analysis of thick rectangular
  plates.
\newblock {\em International Journal of Solids and Structures}, 30:3357--3379,
  19993.

\bibitem{liewzhao2011}
KM~Liew, Xin Zhao, and AJM.
\newblock A review of meshless methods for laminated and functionally graded
  plates and shells.
\newblock {\em Composite Structures}, 93:2031--2041, 2011.

\bibitem{liuwang2010}
DY~Liu, CY~Wang, and WQ~Chen.
\newblock {Free vibration of FGM plates with in-plane material inhomogeneity}.
\newblock {\em Composite Structures}, 92:1047--1051, 2010.

\bibitem{luchen2008}
CF~L\"u, WQ~Chen, RQ~Xu, and CW~Lim.
\newblock Semi analytical elasticity solutions for bi-directional functionally
  graded beams.
\newblock {\em International Journal of Solids and Structures}, 45:258--275,
  2008.

\bibitem{mansong2012}
H~Man, C~Song, , and W~Gao.
\newblock {A unified 3D-based technique for plate bending analysis using scaled
  boundary finite element method}.
\newblock {\em International Journal for Numerical Methods in Engineering},
  91:491--515, 2012.

\bibitem{mansong2013}
H~Man, C~Song, T~Xiang, W~Gao, and F~Tin-Loi.
\newblock High-order plate bending analysis based on the scaled boundary finite
  element method.
\newblock {\em International Journal for Numerical Methods in Engineering},
  95:331--360, 2013.

\bibitem{nemat-alla2003}
A~Nemat-Alla.
\newblock Reduction of thermal stresses by developing two dimensional
  functionally graded materials.
\newblock {\em International Journal of Solids and Structures}, 40:7339--7356,
  2003.

\bibitem{nguyen-xuantran2012}
H~Nguyen-Xuan, Loc~V Tran, Chien~H Thai, and T~Nguyen-Thoi.
\newblock Analysis of functionally graded plates by an efficient finite element
  method with node-based strain smoothing.
\newblock {\em Thin-Walled Structures}, 54:1--18, 2012.

\bibitem{qianbatra2005}
LF~Qian and RC~Batra.
\newblock Design of bidirectional functionally graded plate for optimal natural
  frequencies.
\newblock {\em Journal of Sound and Vibration}, 280:415--424, 2005.

\bibitem{qianching2004}
LF~Qian and HK~Ching.
\newblock Static and dynamic analysis of {2D} functionally graded elasticity by
  using meshless local {Petrov-Galerkin} method.
\newblock {\em J Chinese Inst Eng}, 27:491--503, 2004.

\bibitem{reddy1984}
JN~Reddy.
\newblock A simple higher order theory for laminated composite plates.
\newblock {\em Journal of Applied Mechanics - Transactions of ASME},
  51:745--752, 1984.

\bibitem{reddy2000}
JN~Reddy.
\newblock Analysis of functionally graded plates.
\newblock {\em International Journal for Numerical Methods in Engineering},
  47:663--684, 2000.

\bibitem{reddychin1998}
JN~Reddy and CD~Chin.
\newblock Thermomechanical analysis of functionally graded cylinders and
  plates.
\newblock {\em Journal of Thermal Stresses}, 21:593--629, 1998.

\bibitem{somashekarprathap1987}
BR~Somashekar, G~Prathap, and C~Ramesh Babu.
\newblock A field-consistent four-noded laminated anisotropic plate/shell
  element.
\newblock {\em Computers \& Structures}, 25:345--353, 1987.

\bibitem{songwolf1997a}
C~Song and JP~Wolf.
\newblock The scaled boundary finite-element methods-alias consistent
  infinitesimal finite-element cell methods for elastodynamics.
\newblock {\em Computer Methods in Applied Mechanics and Engineering},
  147:329--355, 1997.

\bibitem{thainguyen-xuan2012}
CH~Thai, H~Nguyen-Xuan, N~Nguyen-Thanh, T.-H Le, T~Nguyen-Thoi, and T~Rabczuk.
\newblock Static, free vibration and buckling analysis of laminated composite
  reissner-mindlin plates using nurbs-based isogeometric approach.
\newblock {\em International Journal for Numerical Methods in Engineering},
  91:571--603, 2012.

\bibitem{tranferreira2013}
LV~Tran, AJM Ferreira, and H~Nguyen-Xuan.
\newblock Isogeometric analysis of functionally graded plates using
  higher-order shear deformation theory.
\newblock {\em Composites Part B: Engineering}, 51:368--383, 2013.

\bibitem{uymazaydogdu2012}
B~Uymaz, M~Aydogdu, and S~Filiz.
\newblock {Vibration analysis of FGM plates with in-plane material
  inhomogeneity by Ritz method}.
\newblock {\em Composite Structures}, 94:1398--1405, 2012.

\bibitem{valizadehnatarajan2013}
N~Valizadeh, S~Natarajan, OA~Gonzalez-Estrada, T~Rabczuk, TQ~Bui, and SPA
  Bordas.
\newblock Nurbs-based finite element analysis of functionally graded plates:
  static bending, vibration, buckling and flutter.
\newblock {\em Composite Structures}, 99:309--326, 2013.

\bibitem{wang1997}
S~Wang.
\newblock Buckling analysis of skew fibre-reinforced composite laminated plates
  based on first-order shear deformation theory.
\newblock {\em Composite Structures}, 37:5--19, 1997.

\bibitem{zinekiewicztaylor2000}
OC~Zinekiewicz and RL~Taylor.
\newblock {\em {The finite element method}}.
\newblock Butterworth Heinemann, Oxford, 2000.

\end{thebibliography}

\end{document}